\documentclass[12pt]{article}
  
\newtheorem{theorem}{Theorem}
 
\newcommand{\p}{\partial}
\newcommand{\om}{\omega}
\newcommand{\e}{\epsilon}

\begin{document}

\title{Time-Periodic Quasigeostrophic Motion 
	under Dissipation and Forcing
\footnote{This work was supported by the National Science
		Foundation Grant DMS-9704345.}  }

\author{Jinqiao Duan \\
Department of Mathematical Sciences,\\
Clemson University, Clemson, South Carolina 29634, USA.\\
E-mail: duan@math.clemson.edu, Fax: (864)656-5230.} 

\date{  }

\maketitle

\begin{abstract}

The quasigeostrophic equation is a prototypical geophysical fluid model.
In this paper, we consider time-periodic motions of this model under
dissipation and time-dependent wind forcing. 
We show that when the wind forcing is time-periodic and
the spatial square-integral of the wind forcing is bounded
in time, the full nonlinear quasigeostrophic model  has
time-periodic motions, 
under some conditions on $\beta$ parameter, Ekman number, viscousity
and the domain size.

\bigskip


\end{abstract}

\bigskip
\bigskip

{\bf Short running title:}

	Time-Periodic Quasigeostrophic Motion

\newpage

\section{Introduction}

The two dimensional barotropic quasigeostrophic (QG) equation has been
derived as an approximation of the rotating shallow water equations
by the conventional asymptotic expansion in small Rossby number
(\cite{Pedlosky}). The lowest order approximation gives the
barotropic QG equation, which is also the conservation law for the
zero-th order potential vorticity. Warn et al. \cite{Warn} and Vallis
\cite{Vallis} emphasize that this asymptotic expansion is generally
secular for all but the simplest flows and propose a modified
asymptotic method, which involves expanding only the fast modes. The
barotropic QG equation also emerges at the lowest order in this
modified expansion.

Schochet (\cite{Schochet}) has recently shown that
quasigeostrophy is a valid approximation of  the rotating
shallow  water equations  in the limit of zero Rossby number,
i.e., at asymptotically high rotation rate.  
For related work in the three dimensional baroclinic QG model,
see, for example,  \cite{Beale} ,  
 \cite{Embid_Majda} , and \cite{Bennett} .

It is known that the {\em linearized} two dimensional barotropic QG equation 
without forcing and without dissipation has time-periodic solutions,
and these time-periodic solutions are damped away by Ekman or viscous
dissipation (\cite{Pedlosky}, pages 147 and 236).

In this paper, we show that when the wind forcing is periodic and 
when its spatial square-integral is bounded
in time, the full nonlinear QG equation with Ekman and viscous
dissipation has time-periodic solutions. We use a topological technique from
nonlinear global analysis (\cite{Krasnoselskii}).

	\section{Quasigeostrophy}

The two dimensional barotropic QG equation is (\cite{Pedlosky}, page 234)
\begin{eqnarray}
  \Delta \psi_t + J(\psi, \Delta \psi) + \beta \psi_x
	=  \nu \Delta^2 \psi - r \Delta \psi + f(x,y,t),
\label{qg}
\end{eqnarray}
where $\psi(x,y,t)$ is the stream function, $\beta > 0$ is the
meridional gradient of the Coriolis parameter,  $\nu >0$ is the viscous dissipation constant
and $r>0$ is the Ekman dissipation constant and $f(x,y,t)$ is the wind forcing.
Moreover, $\Delta = \p_{xx}+\p_{yy} $ is the Laplacian operator in the
$x,y$ plane and $J(f, g) =f_xg_y -f_yg_x$ is the Jacobi operator. 
This equation is written in the often-studied case of infinite
Rossby deformation radius and with flat bottom, applicable for
planetary-scale solutions. The situation with infinite Rossby
deformation radius is equivalent to the rigid-lid approximation 
(\cite{Benoit}, p223).
 
In this paper, we assume that the wind forcing $f(x,y,t)$ is 
periodic in time with period $T>0$.
 
Introducing the relative vorticity 
$\om (x,y,t) = \Delta \psi(x, y,t)$, the equation (\ref{qg}) can
be written as
\begin{equation}
       \om_t + J(\psi, \om ) + \beta \psi_x
        =\nu \Delta \om - r \om  +  f(x,y,t) \; ,
 \label{eqn}
\end{equation}
where $(x, y) \in D$, an arbitrary bounded planar domain
with piecewise smooth boundary.  This equation is supplemented by zero
Dirichlet boundary conditions for both $\psi$
and $\om = \Delta \psi$, together with an appropriate initial condition,
i.e., we require
\begin{eqnarray}
 \psi(x,y,t) & = & 0 \quad \mbox{on} \; \partial D \; , \label{BC1} \\
 \om (x,y,t) & = & 0 \quad \mbox{on} \; \partial D \; , \label{BC2} \\
 \om (x,y,0) & = & \om_0(x,y)                    \; . \label{IC}
\end{eqnarray} 
These boundary conditions have been used in
analytical and numerical study of this model
in, e.g., \cite{Cessi}, \cite{Pierrehumbert},\cite{Barcilon}.

We note that the Poincar\'e inequality (\cite{Gilbarg-Trudinger})
\begin{eqnarray}
\int_D g^2(x,y) dxdy \leq  \frac{|D|}{\pi} \int_D |\nabla g|^2 dxdy
\end{eqnarray}
holds with these boundary
conditions, where $|D|$ is the area of the domain $D$. 
The global well-posedness (smooth solutions) of this dissipative model 
can be obtained similarly as  in, for example, 
\cite{Barcilon}, \cite{Wolansky} and \cite{Wu}.

It is well-known that the {\em linearized} two dimensional 
barotropic QG equation 
without forcing and without dissipation  
\begin{equation}
   \Delta \psi_t  + \beta \psi_x =0, \qquad x, y \in [0,1],
\label{linQG}
\end{equation}
has time-periodic solutions (\cite{Pedlosky}, p.146-149). 
For example, assume the basin is the unit square on the $\beta-$plane,
with boundary condition $\psi = 0$.
The boundary is a streamline in this case.
By separating variables, one can find that the equation (\ref{linQG})  
has time-periodic solutions (\cite{Pedlosky}, p.146-149)
\begin{eqnarray}
\psi_{mn}(x,y,t) = cos(\frac{\beta x}{2 \sigma_{mn}} +\sigma_{mn}t)
			sin(m\pi x) sin(n\pi y),
		m, n =1, 2, 3, \cdots,
\end{eqnarray}
which are basin-scale traveling-wave oscillations with dispersion
relation,
\begin{eqnarray}
\sigma_{mn} = \frac{-\beta}{2\pi \sqrt{ m^2+n^2 }}\ .
\label{linear_waves}
\end{eqnarray}
These are basin-scale normal modes (planetary waves) for QG in the
rigid-lid approximation. Each mode consists of a carrier wave
$cos(\frac{\beta x}{2 \sigma_{mn}} +\sigma_{mn}t)$ moving to the left
(westward) and modulated by an envelope of sine functions which
maintain the boundary conditions.  These time-periodic linear solutions
are damped away by Ekman or viscous dissipation (\cite{Pedlosky}, p. 236).

We address the question of whether there are any  
basin-scale time-periodic solutions in the full nonlinear
dissipative QG dynamics  (\ref{qg})  with time-periodic wind forcing.

	\section{Dissipativity}

In the following we use the abbreviations $L^2=L^2(D)$, $H^1_0 =H^1_0(D)$,
$\ldots$ for the standard Sobolev spaces.  Furthermore, let $<\cdot, \cdot>$,
$\| \cdot \|$ denote the standard scalar product and norm in $L^2 $,
respectively.  We need the following inequalities (\cite{Temam}).  

{\em Young inequality}:
\begin{equation}
AB \leq \frac{\e}2 A^2 + \frac{1}{2 \e}B^2,
\end{equation}
where $A,B$ are non-negative real numbers and $\e >0$.

{\em Gronwall inequality}:
If an integrable function $y(t)$ satisfies that
\begin{equation}
\frac{dy}{dt} \leq A y +B, \;\;\;\;  t\geq 0,
\end{equation}
for some constants $A$, $B$ with $A \neq 0$, then
\begin{equation}
y(t)\leq [ y(0) + \frac{B}{A} ] e^{At} - \frac{B}{A}, \;\;\; t>0.
\end{equation}

We now show that the system (\ref{eqn}), under boundary conditions
(\ref{BC1}), (\ref{BC2}), is a dissipative system in the sense 
(\cite{Krasnoselskii} or \cite{Hale}) that
all solutions $\om(x,y,t)$ approach a bounded set in $L^2(D)$ as
time goes to infinity. A $T-$time-periodic dissipative
system in a Banach space has at least one $T-$time-periodic
solution. This result follows from a Leray-Schauder topological degree
argument and the Browder's principle (\cite{Krasnoselskii}, p.235).    

Multipling (\ref{eqn}) by $\om$ and integrating over $D$, we get

\begin{eqnarray}
&   &  \frac12 \frac{d}{dt}\|\om\|^2 +\int_D J(\psi, \om)\om dxdy
 +\beta \int_D \psi_x \om dxdy   \nonumber   \\
& = &  -\nu \int_D |\nabla \om|^2 -r \| \om \|^2 +\int_D f(x,y,t)\om dxdy.
 	\label{estimate1}
\end{eqnarray} 
Note that 
\begin{eqnarray}
\int_D J(\psi, \om)\om dxdy =0,
	\label{estimate0}
\end{eqnarray}
via integration by parts;
see also \cite{Barcilon} or \cite{Holm}. Moreover,
using the Young and Poincar\'e inequalities, we get

\begin{eqnarray} 
 | \beta \int_D \psi_x \om dxdy |
& \leq &  \frac12 \beta (\int_D \psi_x^2 + \int_D \om^2 )   \nonumber   \\
& \leq &  \frac12 \beta (\frac{|D|}{\pi}\int_D \om^2 dxdy + \int_D \om^2 dxdy)  \nonumber \\
&   =  &  \frac12 \beta (\frac{|D|}{\pi} +1) \int_D \om^2 dxdy,  
	\label{estimate2}  \\
-\nu \int_D |\nabla \om|^2  dxdy
& \leq & -\frac{\pi \nu}{|D|} \int_D \om^2 dxdy. 
	\label{estimate3}   
\end{eqnarray} 

We further assume that the square-integral of the wind forcing $f(x,y,t)$   
with respect to $x, y$ is bounded in time.
\begin{eqnarray} 
\int_D f(x,y,t)\om dxdy 
 \leq   \frac1{\e} \int_D f^2(x,y,t) dxdy + \e \int_D \om^2 dxdy
 \leq   M + \e \int_D \om^2 dxdy, 
 	\label{estimate4}
\end{eqnarray} 
where $\e > 0$ is to be determined, and $ M>0 $ is a time-independent constant 
since we have assume that the
square-integral of $f(x,y,t)$ with respect to $x, y$ is bounded in time.

Putting (\ref{estimate0}), (\ref{estimate2}), (\ref{estimate3}), 
(\ref{estimate4}) into (\ref{estimate1}), we obtain
\begin{eqnarray}
 \frac12 \frac{d}{dt}\|\om\|^2 +  \alpha \int_D \om^2 dxdy  \leq   M,
\end{eqnarray} 
where
$$
\alpha = [r+\frac{\pi \nu}{|D|}-\frac12 \beta (\frac{|D|}{\pi} +1) -\e ].
$$
Assume that
\begin{eqnarray}
    r+\frac{\pi \nu}{|D|} > \frac12 \beta (\frac{|D|}{\pi} +1).
\end{eqnarray}
We can then take $\e>0$ small enough such that $\alpha > 0$.
Thus, by the Gronwall inequality, we have
\begin{eqnarray}
 \|\om\|^2 \leq (\|\om_0\|^2 -\frac{M}{\alpha})e^{-2\alpha t} 
 	+ \frac{M}{\alpha}.
\end{eqnarray} 
Hence all solutions $\om$ enter a bounded
set $\{ \om:  \;\;   \|\om\|   \leq   \sqrt{\frac{M}{\alpha}} \}$
as time goes to infinity. The system (\ref{eqn}) is  therefore a
dissipative system and hence has at least one $T-$time-periodic
solution.

We then have the following result.

\begin{theorem}
Assume that

(i) the wind forcing $f(x,y,t)$ is time-periodic with period $T>0$,
and its square-integral with respect to $x, y$ is bounded in time; and

(ii) $r+\frac{\pi \nu}{|D|} > \frac12 \beta (\frac{|D|}{\pi} +1)$,
where $\beta > 0$ is the meridional gradient of the Coriolis parameter,  
$\nu >0$ is the viscous dissipation constant, 
$r>0$ is the Ekman dissipation constant, and $|D|$ is the area of the 
bounded domain $D$.

Then the dissipative quasigeostrophic model
\begin{eqnarray}
  \Delta \psi_t + J(\psi, \Delta \psi) + \beta \psi_x
       & = & \nu \Delta^2 \psi - r \Delta \psi + f(x,y,t) \; ,  \\
  \om   & = &\Delta \psi \; ,  \\
  \psi(x,y,t) & = & 0 \quad \mbox{on} \; \partial D \; ,  \\
  \om (x,y,t) & = & 0 \quad \mbox{on} \; \partial D \; ,   \\
  \om (x,y,t+T) & = & \om (x,y,t) \; ,
\end{eqnarray} 
has at least one  time-periodic solution with period $T>0$.
\end{theorem}

We   remark that it is generally difficult to show
existence of time-periodic motions for a spatially extended
evolution system. Our result provides such a proof of existence, 
under some conditions on $\beta$ parameter, Ekman number, viscousity
and the domain size, for a prototypical geophysical fluid model.

\bigskip

{\bf Acknowledgement.} This research was performed while the 
author was visiting the Isaac Newton Institute for 
Mathematical Sciences at Cambridge University.

\end{document}